\documentclass{article}
\usepackage{amsmath,amssymb,amsfonts}

\newtheorem{theorem}{Theorem}
\newtheorem{proposition}{Proposition}
\newtheorem{corollary}{Corollary}

\def\R{{\mathbb R}}

\def\Re{{\mathrm{Re}\, }}
\def\Im{{\mathrm{Im}\, }}
\def\nil{{\mathrm{Nil}\, }}
\def\sol{{\mathrm{Sol}\, }}
\def\sl{\widetilde{SL}_2}
\def\sll{{\mathrm{SL}\, }}

\begin{document}

\title{Surfaces in three-dimensional Lie groups
\thanks{The work is supported by RFBR (grant no. 03-01-00403) and
by the programm of fundamental researches of RAS ``Mathematical methods in nonlinear dynamics''.}
}
\author{Dmitry A. Berdinsky
\thanks{Department of Mechanics and Mathematics, Novosibirsk State University,
630090 Novosibirsk, Russia; e-mail: berdinsky@ngs.ru
.}
\and
Iskander A. Taimanov
\thanks{Institute of Mathematics, 630090 Novosibirsk, Russia;
e-mail: taimanov@math.nsc.ru.}
}
\date{}
\maketitle

\section{Introduction}
\label{sec1}

In the present paper we extend the methods of the Weierstrass (or spinor) representation of surfaces in $\R^3$
\cite{T1,T2} and $SU(2)=S^3$ \cite{T3} for surfaces in the three-dimensional Lie groups $\nil, \sl$, and $\sol$
endowed with the so-called Thurston's geometries \cite{Scott}.

The main feature of this approach is that
the geometry of a surface is related to the spectral properties of the corresponding Dirac operator.
Therewith this approach reveals some
unknown  before geometric meanings of the Willmore functional and the Willmore conjecture
which states that for tori the Willmore functional attains its minimum on
the Clifford torus.

A function $\psi$ generates a surface in $\R^3$ via the Weierstrass formulas
if and only if it meets some equation of the Dirac type where the Dirac operator in general has
two potentials $U$ and $V$ which coincide for the case of surfaces in $\R^3$.
Given a surface $M$ in $\R^3$, the integral
$$
E(M) = \int_M UV \frac{i dz \wedge d\bar{z}}{2}
$$
equals to
$\frac{1}{4}\int_M H^2 d\mu$
where $H$ is the mean curvature and $d\mu$ is the induced area form on $M$,
i.e. $E(M)$ is equal up to a multiple to
the Willmore functional
$$
{\cal W}(M) = \int_M H^2 d\mu
$$
(see \cite{T1}).

The same was established for surfaces in $SU(2)$ in \cite{T3} although the Weierstrass formulas have to be replaced
by an analogous construction valid for noncommutative Lie groups.

Looking for a physical confirmation of
the Willmore conjecture we notice that for the Clifford torus the spectral curve
which was defined for general tori in \cite{T2} has zero geometric genus. Moreover the deviation of
the spectral curve of a torus from the flat curve which is
the spectral curve of the Dirac operator with vanishing potentials
$U=V=0$ is measured by $E(M)$. This led us to an approach for proving the Willmore conjecture
by proving that in each conformal class
the minimum of the Willmore functional is attained on a torus with the minimal geometric genus
(this part some time ago was announced to be proved by M.U. Schmidt (see math.DG/0203224 in http://arxiv.org))
as the first step and by checking the conjecture for tori with the minimal geometric genus in their  conformal classes
as the final step.

Recently M. Haskins proposed to treat the genus of the spectral curve as a measure of geometric complexity
for some other variational problems of geometry \cite{Haskins}.

In the present paper we derive the analogous formulas for the functional $E(M)$
for surfaces in $\nil,\sl$, and $\sol$. It appears it has very nonexpectable geometric properties.
We call it the spinor energy or just the energy of a surface.

Until now we even do not know is it bounded from below or not.
For surfaces in $\sol$ we do not prove that it is real-valued.
However it measures the deviation of the spectral curve of a torus from the flat curve and we think
that that is enough for confirming its geometric importance.

The problems of finding the analogs of the Willmore conjecture for such functionals and describing their extremals
are very interesting.

In the present paper we also derive the equations for $\psi$ corresponding to minimal surfaces in Lie groups
and by our means obtain another proof to Abresch's result that for constant mean curvature surfaces in
$\nil$ and $\sl$ certain quadratic differentials are holomorphic \cite{Abresch}.
It would be very intriguing to relate these conditions
to integrable systems as it was done for constant mean curvature surfaces in $\R^3$ and $S^3$.

We note that the approach by Abresch differs from our and bases on a representation of $\nil$ and $\sol$ as
line bundles of constant curvature over surfaces of constant curvature (see also the recent paper \cite{Daniel}
where this property is used for studying surfaces $\nil$ and $\sol$). It looks that the failure of such a
representation for $\sol$ explains why our approach does not lead to such satisfactory understanding of surfaces as
in  the cases of $\nil$ and $\sl$.

We notice that the study of integrable surfaces in arbitrary Lie groups by methods of integrable systems
was first attempted in \cite{FG}. We also hope that our approach will be helpful in studying
global properties of minimal and constant mean curvature surfaces in Lie groups in the spirit of \cite{AR,FMP}.

The results of this paper were partially exposed on the conference on the surface theory in Benediktbeuern (January 2005).

We thank U. Abresch for many helpful conversations and an explanation of his result on holomorphic quadratic
differentials for constant mean curvature surfaces in $\nil$ and $\sol$.

\section{Preliminary facts}
\label{sec2}

\subsection{Left invariant metrics on Lie groups}
\label{ssec2.1}

Recall that a metric $\langle \xi,\eta\rangle$ on  Lie group $G$
is called left invariant if it is invariant with respect to left
translations:
$$
L_g: G \stackrel{\times g}{\longrightarrow} G \ \ \ : \ \ \ h \to
gh, \ h \in G,
$$
i.e. for all vectors $\xi$ and $\eta$ tangent to $G$ at $h$,
the inner product of their translations by $g$,
and the inner product of  $\xi$ and $\eta$ coincide:
$$
\langle \xi, \eta \rangle = \langle L_g^\ast \xi, L_g ^\ast \eta
\rangle, \ \ \ \xi, \eta \in T_hG, \ g,h \in G.
$$
Clearly every left invariant metric is determined by the inner
product of vectors tangent at the unit $1 \in G$ of the group $G$,
i.e. by the inner product on the Lie algebra ${\cal G}$ of $G$.

For calculating the Levi-Civita connection we apply the general
tetrad formalism known from mathematical physics which is as
follows.

Given vector fields $e_1,\dots,e_n$ on an $n$-dimensional manifold
$M$ such that at every point the corresponding vectors form an
orthonormal base:
$$
\langle e_i,e_j \rangle = \delta_{ij}, \ \ \ i,j =1, \dots,n,
$$
the Levi-Civita connection is given by the following formulas:
$$
\Gamma^i_{jk} = \frac{1}{2} \left(c^i_{kj} + c^j_{ik} +
c^k_{ij}\right),
$$
$$
\nabla_{e_k} e_j = \Gamma^i_{jk} e_i,
$$
where
$$
[e_i, e_j] = c^k_{ij} e_k.
$$

In our case let us take an orthonormal base $\xi_1,\dots,\xi_n$
for the tangent space at the unit of the group $G$ and extend it
to vector fields on the whole group by left translations:
$$
e_i(g) = L_g^\ast \xi_i, \ \ \ i=1,\dots,n, \ g \in G.
$$
Then $c^i_{jk}$ are just the structure constants of the Lie
algebra ${\cal G}$ of the group $G$ and putting for simplicity
$$
\alpha_{ijk} = \langle [e_i,e_j],e_k \rangle = c^k_{ij}
$$
we derive
\begin{equation}
\label{connection}
\nabla_{e_k}e_j = \frac{1}{2} \sum_i
\left(\alpha_{kji} + \alpha_{ikj} + \alpha_{ijk} \right) e_i.
\end{equation}
In particular,
$\alpha_{ijk}$ is a skew-symmetric tensor for
a compact Lie group $G$ with the Killing
metric and the formula \eqref{connection}
reduces to
$$
\nabla_X Y = \frac{1}{2}[X,Y]
$$
where $X$ and $Y$ are left invariant vector fields.
For details we refer to \cite{Milnor}.

\subsection{The derivational equations}
\label{ssec2.2}

Let $\Sigma$ be a surface immersed into $G$ and let
$$
f: \Sigma \to G
$$
be the immersion. Choose a conformal parameter $z=x+iy$ on $\Sigma$
(or, more precisely, in a domain of $\Sigma$)
and denote by ${\bf I} = e^{2\alpha} dz d\bar{z}$ the induced metric.

Let us consider the pullback of $TG$ to a  ${\cal G}$-bundle
over $\Sigma$:
${\cal G} \to E = f^{-1}(TG) \stackrel{\pi}{\to} \Sigma$
and the differential
$$
d_{\cal A}:\Omega^1(\Sigma;E) \to \Omega^2(\Sigma;E),
$$
which acts on $E$-valued $1$-forms as follows. Let us write down
a form $\omega$ as
$$
\omega = u dz + u^\ast d\bar{z}.
$$
Then
$$
d_{\cal A} \omega = d'_{\cal A} \omega + d''_{\cal A} \omega
$$
where
$$
d'_{\cal A} \omega = -\nabla_{\bar{\partial}f} u dz \wedge d\bar{z},
\ \ \
d''_{\cal A} \omega = \nabla_{\partial f}u^\ast dz \wedge d\bar{z}.
$$

By straightforward computations we obtain the first derivational equation
\begin{equation}
\label{d1}
d_{\cal A} (df) = 0.
\end{equation}

By the definition of the tension vector $\tau(f)$, we have
$$
d_{\cal A} (\ast df) = f\cdot (e^{2\alpha} \tau(f)) dx \wedge dy =
\frac{i}{2} f\cdot(e^{2\alpha} \tau(f)) dz \wedge d\bar{z}
$$
where $f\cdot\tau(f) = 2HN$, $N$ is the normal vector and $H$ is the mean curvature.
Hence we obtain the second derivational equation:
\begin{equation}
\label{d2}
d_{\cal A} (\ast df) = i e^{2\alpha} H N dz \wedge d\bar{z}.
\end{equation}

For $S=SU(2)$ we expose this scheme in \cite{T3}
(see also the derivation of the harmonicity equation for surfaces in Lie groups
in \cite{Hitchin}).

Since the metric is left invariant we write down these equations in terms of
$$
\Psi = f^{-1}\partial f, \ \ \
\Psi^\ast = f^{-1} \bar{\partial} f
$$
as
\begin{equation}
\label{h1}
\partial\Psi^\ast - \bar{\partial}\Psi + \nabla_{\Psi}\Psi^\ast -
\nabla_{\Psi^\ast}\Psi = 0,
\end{equation}
\begin{equation}
\label{h2}
\partial\Psi^\ast + \bar{\partial}\Psi + \nabla_{\Psi}\Psi^\ast +
\nabla_{\Psi^\ast}\Psi = e^{2\alpha} H f^{-1}(N).
\end{equation}
The equation \eqref{h1} is equivalent to \eqref{d1} and
the equation \eqref{h2} is equivalent to \eqref{d2}.

In the sequel we assume that the Lie group $G$ is three-dimensional and
we choose an orthonormal basis $e_1,e_2,e_3$ for the inner product on the
Lie algebra ${\cal G}$ of the Lie group $G$.

Let us expand $\Psi$ and $\Psi^\ast$ in this basis as
$$
\Psi = \sum_{k=1}^3 Z_k e_k, \ \ \
\Psi^\ast = \sum_{k=1}^3 \bar{Z}_k e_k
$$
and rewrite \eqref{h1} and \eqref{h2} in terms of
$Z_k$, $k=1,2,3$, as follows
\begin{equation}
\label{z1}
\sum_j (\partial \bar{Z}_j - \bar{\partial} Z_j) e_j +
\sum_{j,k} ( Z_j \bar{Z}_k - \bar{Z}_j Z_k) \nabla_{e_j}e_k = 0,
\end{equation}
\begin{equation}
\label{z2}
\begin{split}
\sum_j (\partial \bar{Z}_j + \bar{\partial} Z_j) e_j +
\sum_{j,k} (Z_j \bar{Z}_k + \bar{Z}_j Z_k) \nabla_{e_j}e_k = \\
2iH
\left[
(\bar{Z}_2 Z_3 - Z_2 \bar{Z}_3) e_1 +
(\bar{Z}_3 Z_1 - Z_3 \bar{Z}_1) e_2 +
(\bar{Z}_1 Z_2 - Z_1 \bar{Z}_2) e_3
\right].
\end{split}
\end{equation}
Here we assumed that the basis $\{e_1,e_2,e_3\}$ is positively
oriented and therefore
$$
f^{-1}(N) = 2 i  e^{-2\alpha}
\left[
(\bar{Z}_2 Z_3 - Z_2 \bar{Z}_3) e_1 +
(\bar{Z}_3 Z_1 - Z_3 \bar{Z}_1) e_2 +
(\bar{Z}_1 Z_2 - Z_1 \bar{Z}_2) e_3
\right]
$$
(for $G=SU(2)$ with the Killing metric this formula takes the form
$f^{-1}(N)= 2ie^{-2\alpha} [\Psi^\ast,\Psi]$. \footnote{In
\cite{T3} the term of $2$ in the right hand side of the same
formula for $N$ (the formula (47)) was skipped however it was
counted in all computations and this misprint does not affect any result.}
Since the parameter $z$ is conformal we have
$$
\langle \Psi,\Psi \rangle = \langle \Psi^\ast,\Psi^\ast\rangle = 0,
\ \ \
\langle \Psi, \Psi^\ast\rangle = \frac{1}{2} e^{2\alpha}
$$
which is rewritten as
$$
Z_1^2 + Z_2^2 + Z_3^2 = 0, \ \ \
|Z_1|^2+|Z_2|^2+|Z_3|^2 = \frac{1}{2} e^{2\alpha}.
$$
The first equality implies that the vector $Z$ can be
parameterized in the form
\begin{equation}
\label{spinor}
Z_1 = \frac{i}{2} ( \bar{\psi}_2^2 + \psi_1^2), \ \ \
Z_2 = \frac{1}{2} ( \bar{\psi}_2^2 - \psi_1^2), \ \ \
Z_3 = \psi_1 \bar{\psi}_2.
\end{equation}

\subsection{The Dirac operator and the energy of a surface}
\label{ssec2.3}

The Weierstrass representation is derived by inserting $\psi$ into
\eqref{z1} and \eqref{z2} (see  \cite{T1,T3} for such representations
of surfaces in $\R^3$ and $S^3$).

In this event the derivational equations takes the form of the Dirac equation
\begin{equation}
\label{dirac}
{\cal D} \psi =
\left[\left(
\begin{array}{cc}
0 & \partial \\
-\bar{\partial} & 0
\end{array}
\right)
 +
\left(
\begin{array}{cc}
U & 0 \\
0 & V
\end{array}
\right)
\right]
\psi = 0
\end{equation}
which is satisfied by
$$
\psi = \left(\begin{array}{c} \psi_1 \\ \psi_2 \end{array}\right).
$$
The first fundamental form is equal to
$$
{\bf {\mathrm I}} = e^{2\alpha} dz d\bar{z} =
(|\psi_1|^2+|\psi_2|^2)^2 dz d\bar{z}
$$
and the left translation of the normal vector is
\begin{equation}
\label{normal}
\begin{split}
f^{-1}(N) =
e^{-\alpha}\left[ i(\psi_1 \psi_2 - \bar{\psi}_1\bar{\psi}_2) e_1 -
\right.
\\
\left.
(\psi_1 \psi_2 + \bar{\psi}_1\bar{\psi}_2) e_2 +
(|\psi_2|^2-|\psi_1|^2) e_3)\right].
\end{split}
\end{equation}

Since $\Psi dz = (\sum Z_k e_k) dz$ is a correctly defined
$1$-form , we conclude that
the $1$-forms
$$
\psi_1^2 dz, \ \ \bar{\psi}_2^2 dz, \ \ \psi_1 \bar{\psi}_2 dz
$$
are globally defined on the whole surface. By the Dirac equation,
this implies that
the potentials $U$ and $V$ are the squares of $dz d\bar{z}$ forms
and therefore
$$
UV dz \wedge d\bar{z}
$$
is a well-defined $2$-form on the surface (see \cite{T1} for details).

For a surface $M$ immersed into $\R^3$ the potentials of the Dirac
operator take real values and coincide: $U=V$. For such a surface
the following equality holds:
$$
E(M) = \int_M UV \frac{i dz \wedge d\bar{z}}{2} = \frac{1}{4}
\int_M H^2 d\mu
$$
where $d\mu = e^{2\alpha} dx\wedge dy$ is the induced measure on the surface and
$H$ is the mean curvature. We recall that the
Willmore functional
equals
$$
{\cal W}(M) = \int_M H^2 d\mu.
$$
This observation from \cite{T1}
was the starting point for treating the spectral quantities of
this operator $\cal D$ as the geometric quantities and for
physical explanation of the Willmore conjecture by using the
spectral curves of immersed tori (see \cite{T1,T3}).

Generically $U$ and $V$ do not coincide necessarily and we have to
consider the energy of a closed surface as
\begin{equation}
\label{energy}
E(M) = \int_M UV dx \wedge dy = \int_M UV  \frac{i dz \wedge d\bar{z}}{2}.
\end{equation}

It is even not clear that it always be real-valued. We shall check
that for some special cases. For that we shall use the following simple
proposition which follows from the Dirac equation.
Indeed we have
$$
\bar{\partial} \psi_1 = (\Re V + i \Im V) \psi_2,\ \ \
\bar{\partial}\bar{\psi}_2 = (-\Re U +i \Im U)\psi_1
$$
which implies

\begin{proposition}
\label{identity}
Given $\psi$ meeting the Dirac equation \eqref{dirac},
the following identity holds
\begin{equation}
\label{eq-identity}
\bar{\partial}(\psi_1 \bar{\psi}_2) =
\left( -\Re U|\psi_1|^2 + \Re V |\psi_2^2|\right)
+i
\left(\Im U|\psi_1|^2 + \Im V|\psi_2|^2 \right).
\end{equation}
\end{proposition}

\subsection{Thurston's geometries on Lie groups $\nil, \sol$ and $\sl$}
\label{ssec2.4}

By Thurston's theorem all three-dimensional maximal simply connected
geometries $(X,\mathrm{Isom}\, X)$ admitting compact quotient geometries
belong to the following list:

1) three geometries with constant sectional curvature:
$X = \R^3, S^3$, or $H^3$;

2) a pair of product geometries: $X = S^2 \times \R$ or $H^2 \times \R$;

3) three geometries modeled on the Lie groups $\nil, \sol$, and $\sl$
with certain left invariant metrics.

We refer for details of this classification theorem and for the related
famous Thurston geometrization conjecture for three-manifolds
to \cite{Thurston,Scott}.

In the paper we study surfaces in these geometries modelled on noncompact
Lie groups $\nil, \sol$, and $\sl$. Before going to the surface theory
we briefly recall main facts about these geometries referring for
an advanced exposition to \cite{Scott}.

\subsubsection{The group $\nil$}
This group is formed by
all matrices of the form
$$
\left(
\begin{array}{ccc}
1 & x & z \\
0 & 1 & y \\
0 & 0 & 1
\end{array}
\right), \ \ \ \ x,y,z \in \R,
$$
with the usual multiplication rule and the left invariant metric
$$
ds^2 = dx^2 + dy^2 + (dz - x dy)^2.
$$
The Lie algebra is generated by the elements
$$
e_1 =
\left(
\begin{array}{ccc}
0 & 1 & 0 \\
0 & 0 & 0 \\
0 & 0 & 0
\end{array}
\right), \ \
e_2 =
\left(
\begin{array}{ccc}
0 & 0 & 0 \\
0 & 0 & 1 \\
0 & 0 & 0
\end{array}
\right), \ \
e_3 =
\left(
\begin{array}{ccc}
0 & 0 & 1 \\
0 & 0 & 0 \\
0 & 0 & 0
\end{array}
\right),
$$
which meet the following commutation relations
$$
[e_1,e_2] = e_3, \ \ [e_1,e_3]=[e_2,e_3]=0.
$$

We derive from \eqref{connection} that the connection is as follows:
\begin{equation}
\label{nil-con}
\begin{split}
\nabla_{e_1}e_2 = - \nabla_{e_2}e_1 = \frac{1}{2}e_3, \ \ \
\nabla_{e_1}e_3 = \nabla_{e_3}e_1 = -\frac{1}{2} e_2, \\
\nabla_{e_2}e_3 = \nabla_{e_3}e_2 = \frac{1}{2}e_1, \ \ \
\nabla_{e_1}e_1 = \nabla_{e_2}e_2 = \nabla_{e_3}e_3 = 0.
\end{split}
\end{equation}

\subsubsection{The group $\sl$}

The group $G = \sl$ is the universal covering of the group $SL(2)$
formed by all real $2 \times 2$-matrices with unit determinant.
The group $PSL(2) = SL(2)/\pm 1$ is the group of orientation
preserving isometries of the hyperbolic plane $H^2$ and it
is diffeomorphic to the unit tangent bundle $UH^2$ of $H^2$.
Moreover the left invariant metric on $\sl$ corresponding to one
of the Thurston geometries (i.e., with the maximal isometry group)
is the pullback under the projection
$$
\sl \to SL(2) \to PSL(2) \approx UH^2
$$
of the metric on $UH^2$ which is as follows.

Let us realize the hyperbolic plane as the disc $|z|<1$ on
the complex plane with the metric
$$
ds^2 = \frac{4 dz d\bar{z}}{(1-|z|^2)^2}
$$
and take on $UH^2$ the metric
$$
dl^2 = \frac{4 (dz d\bar{z} + d\varphi^2)}{(1-|z|^2)^2}
$$
where $\varphi$ is the angle coordinate on the unit circles
in the tangent spaces to $H^2$.
We identify the generators of the Lie algebra of $\sl$ with the
generators of the following one-parametric subgroups
$$
g_x =
\left(
\begin{array}{cc}
\frac{1}{\sqrt{1-x^2}} &
\frac{x}{\sqrt{1-x^2}} \\
\frac{x}{\sqrt{1-x^2}} &
\frac{1}{\sqrt{1-x^2}}
\end{array}
\right), \ \
g_y =
\left(
\begin{array}{cc}
\frac{1}{\sqrt{1-y^2}} &
\frac{iy}{\sqrt{1-y^2}} \\
-\frac{iy}{\sqrt{1-y^2}} &
\frac{1}{\sqrt{1-y^2}}
\end{array}
\right), \ \
$$
$$
g_\varphi =
\left(
\begin{array}{cc}
e^{i\varphi/2} & 0 \\ 0 & e^{-i\varphi/2}
\end{array}
\right)
$$
which act on $H^2$ by fractional linear transformations
$$
z \to \frac{az+b}{cz+d}, \ \ \
\left(\begin{array}{cc} a & b \\ c & d \end{array}\right) \in SL(2).
$$
The generators of these subgroups are
$$
f_1 =
\left(
\begin{array}{cc}
0 & 1 \\ 1 & 0
\end{array}\right), \ \
f_2 =
\left(
\begin{array}{cc}
0 & i \\ -i & 0
\end{array}\right), \ \
f_3 =
\left(
\begin{array}{cc}
\frac{i}{2} & 0 \\ 0 & -\frac{i}{2}
\end{array}\right).
$$
These generators meet the commutation relations
$$
[f_1,f_2] = -4f_3, \ \ \
[f_1,f_3] = -f_2, \ \ \
[f_2,f_3] = f_1.
$$
The desired metric on $\sl$ is induced by the inner product
$$
\langle f_j, f_k \rangle = 4\delta_{jk}, \ \ \ j,k=1,2,3,
$$
on the Lie algebra of $\sl$ and, by \eqref{connection},
this implies
\begin{equation}
\label{sl-con}
\begin{split}
\nabla_{e_1}e_2 = - \nabla_{e_2}e_1 = - e_3, \ \ \
\nabla_{e_1}e_3 = e_2, \ \ \ \nabla_{e_3}e_1 = \frac{3}{2} e_2, \\
\nabla_{e_2}e_3 = -e_1, \nabla_{e_3}e_2 = -\frac{3}{2}e_1, \ \ \
\nabla_{e_1}e_1 = \nabla_{e_2}e_2 = \nabla_{e_3}e_3 = 0,
\end{split}
\end{equation}
where
$$
e_j = \frac{1}{2} f_j, \ \ \
\langle e_j, e_k \rangle = \delta_{jk}, \ \ \
j,k=1,2,3.
$$

\subsubsection{The group $\sol$}

This group consists of
all matrices of the form
$$
\left(
\begin{array}{ccc}
e^{-z} & 0 & x \\
0 & e^z & y \\
0 & 0 & 1
\end{array}
\right), \ \ \ \ x,y,z \in \R,
$$
with the usual multiplication rule and the left invariant metric
$$
ds^2 = e^{2z} dx^2 + e^{-2z} dy^2 + dz^2.
$$
The Lie algebra is generated by the elements
$$
e_1 =
\left(
\begin{array}{ccc}
0 & 0 & 1 \\
0 & 0 & 0 \\
0 & 0 & 0
\end{array}
\right), \ \
e_2 =
\left(
\begin{array}{ccc}
0 & 0 & 0 \\
0 & 0 & 1 \\
0 & 0 & 0
\end{array}
\right), \ \
e_3 =
\left(
\begin{array}{ccc}
-1 & 0 & 0 \\
0 & 1 & 0 \\
0 & 0 & 0
\end{array}
\right),
$$
which satisfy the commutation relations
$$
[e_1,e_2] = 0, \ \ [e_1,e_3]=e_1, \ \ [e_2,e_3]=-e_2.
$$

By \eqref{connection}, we conclude that
\begin{equation}
\label{sol-con}
\begin{split}
\nabla_{e_1}e_2 = \nabla_{e_2}e_1 = 0, \ \ \
\nabla_{e_1}e_3 = e_1, \ \ \ \nabla_{e_3}e_1 = 0, \\
\nabla_{e_2}e_3 = -e_2, \ \ \ \nabla_{e_3}e_2 = 0, \\
\nabla_{e_1}e_1 = - \nabla_{e_2}e_2 = -e_3, \ \ \ \nabla_{e_3}e_3 = 0.
\end{split}
\end{equation}

\subsection{The curvature tensors of three-dimensional Lie groups}
\label{ssec2.5}

By \eqref{nil-con}, \eqref{sl-con}  and \eqref{sol-con} we compute
the curvature tensor
$$
\langle R(X,Y)Z,W \rangle = R_{lkji} W^l Z^k Y^j X^i = \langle
\left(\nabla_Y\nabla_X - \nabla_X \nabla_Y + \nabla_{[X,Y]}\right)
Z, W \rangle
$$
of the groups $\nil,\sl$, and $\sol$. We skip all these simple
computations presenting only the result that we need.

\begin{proposition}
\label{curvatures}
For the groups $\nil,\sl$, and $\sol$
if there are three different indices among $i,j,k,l$ then
$R_{ijkl} = 0$.
The other components of the curvature tensor are as follows
$$
R_{1212} =
\begin{cases}
-\frac{3}{4} & \text{for $\nil$} \\
-4  & \text{for $\sl$} \\
1 & \text{for $\sol$}
\end{cases}
,
\ \ \
R_{1313} = R_{2323} =
\begin{cases}
\frac{1}{4} & \text{for $\nil$} \\
1 & \text{for $\sl$} \\
-1 & \text{for $\sol$}
\end{cases}
.
$$
\end{proposition}

Recall that if $|X|=|Y|=1$ and these vectors are linearly
independent then
$$
K_{XY} = \langle R(X,Y)X,Y \rangle
$$
equals the sectional curvature of the plane spanned by $X$ and
$Y$. This proposition enables us to compute the sectional curvature
of every plane.

\section{The Weierstrass representation for surfaces in Lie groups}
\label{sec3}

\subsection{A construction of a surface from $\psi$}
\label{ssec3.0}

In this section we derive the compatibility conditions for a vector function $\psi$,
i.e. a criterion for $\psi$
to correspond to an immersion of a surface into $\nil, \sl$, or $\sol$.
Given a function $\psi$ meeting these conditions for a Lie group $G$
from this list
we may construct a surface as follows.
\footnote{This generalizes the Weierstrass representation of surfaces in
$\R^3$ \cite{T2} and $SU(2)$ \cite{T3} for surfaces in these groups.}

Let $\psi$ be defined on a surface $M$ with a complex parameter $z$.
Let us pick up a point $P \in M$.

First, we insert $\psi$ into \eqref{spinor} for
the components
$Z_1,Z_2,Z_3$ of $\Psi = \sum_{k=1}^3 Z_k e_k = f^{-1} \partial f$.
Second, we solve the following linear equation in the Lie group $G$:
$$
f_z = f \Psi
$$
with a certain initial data $f(P) = g \in G$.
We thus obtain the desired surface as the mapping
$$
f: M \to G.
$$

It follows from the derivation of the compatibility conditions in \S \ref{ssec2.1} that
any surface $f: M \to G$  is obtained by this construction.
The induced metric takes the form
$$
ds^2 = e^{2\alpha} dz d\bar{z}, \ \ \ e^\alpha  =  (|\psi_1|^2+|\psi_2|^2),
$$
and the induced measure on $M$ is
$$
d \mu = e^{2\alpha} \frac{i dz \wedge d\bar{z}}{2}.
$$
The Hopf quadratic differential equals
$$
\omega = Adz^2, \ \ \ A = \langle \nabla_{f_z} f_z, N \rangle.
$$
Its explicit formula in terms of $\psi$ depends on the Lie group and the Levi-Civita
connection on this group. Indeed, we have
\begin{equation}
\label{hopf}
\begin{split}
A = \langle \Psi_z, N \rangle + \langle \sum_{j,k} Z_j Z_k \nabla_{e_j} e_k, N \rangle =
\\
=
(\bar{\psi}_2 \partial \psi_1 - \psi_1 \partial \bar{\psi}_2) +
\langle \sum_{j,k} Z_j Z_k \nabla_{e_j} e_k, N \rangle.
\end{split}
\end{equation}

The compatibility conditions takes the form of the Dirac equation \eqref{dirac}
and describe $\bar{\partial} \psi_1$ and $\partial \psi_2$:
$$
{\cal D} \psi =
\left[\left(
\begin{array}{cc}
0 & \partial \\
-\bar{\partial} & 0
\end{array}
\right)
 +
\left(
\begin{array}{cc}
U & 0 \\
0 & V
\end{array}
\right)
\right]
\psi = 0.
$$

Another derivatives which are $\partial \psi_1$ and $\bar{\partial}\psi_2$
are derived as follows. We differentiate
$e^\alpha$ and obtain
$$
\alpha_z e^{\alpha} =
\bar{\psi}_1 \partial \psi_1 + \psi_2
\partial \bar{\psi}_2 +( \psi_1 \partial \bar{\psi}_1 +\bar{\psi}_2 \partial \psi_2),
$$
where, by the Dirac equation, the expression in brackets is written as
$$
\psi_1 \partial \bar{\psi}_1 +\bar{\psi}_2 \partial \psi_2 = (\bar{V}-U) \psi_1 \bar{\psi}_2.
$$
Together with the formula for the Hopf differential this gives a system for
$\partial \psi_1$ and $\partial \bar{\psi}_2$:
$$
\left(
\begin{array}{cc}
\bar{\psi}_1 & \psi_1 \\
\bar{\psi}_2 & -\psi_1
\end{array}
\right)
\left(
\begin{array}{c}
\partial \psi_1 \\
\partial \bar{\psi}_2
\end{array}
\right)
=
\left(
\begin{array}{c}
\alpha_z e^\alpha + (U - \bar{V}) \psi_1 \bar{\psi}_2 \\
A - \langle \sum_{j,k} Z_j Z_k \nabla_{e_j} e_k, N \rangle
\end{array}
\right).
$$
Resolving this system we obtain the expressions for
$\partial \psi_1$ and $\bar{\partial} \psi_2$
which together with the Dirac equation give us the complete set of the Weingarten equations
in terms on $\psi$:
\begin{equation}
\label{weingarten}
(\partial - {\cal A})\psi = 0, \ \ \ (\bar{\partial} - {\cal B})\psi = 0.
\end{equation}

The Codazzi equations are now derived from the zero-curvature conditions
$$
\partial \bar{\partial} \psi_k = \bar{\partial} \partial \psi_k, \
\ \ k=1,2.
$$

In the sequel we call $\psi$
the generating spinor of a surface (see \cite{T3} for an explanation of this terminology).

\subsection{The group $\nil$}
\label{ssec3.1}

By \eqref{nil-con}, the derivational equations \eqref{z1} and
\eqref{z2} take the form
\begin{equation}
\label{nil-1}
\begin{aligned}
\partial \bar{Z}_1 - \bar{\partial} Z_1 = 0, \\
\partial \bar{Z}_2 - \bar{\partial} Z_2 = 0, \\
\partial \bar{Z}_3 - \bar{\partial} Z_3 +
(Z_1\bar{Z}_2 - \bar{Z}_1 Z_2) = 0,\\
\partial \bar{Z}_1 + \bar{\partial} Z_1 + (Z_2\bar{Z}_3 + \bar{Z}_2 Z_3) =
2iH(\bar{Z}_2 Z_3 - Z_2 \bar{Z}_3), \\
\partial \bar{Z}_2 + \bar{\partial} Z_2 - (Z_1\bar{Z}_3 + \bar{Z}_1 Z_3)=
2iH(\bar{Z}_3 Z_1 - Z_3 \bar{Z}_1), \\
\partial \bar{Z}_3 + \bar{\partial} Z_3 =
2iH(\bar{Z}_1 Z_2 - Z_1 \bar{Z}_2).
\end{aligned}
\end{equation}
Inserting \eqref{spinor} into these formulas we rewrite
the first pair of equations as
$$
\partial \psi^2_2 + \bar{\partial} \psi_1^2 = 0,
$$
and the pair consisting of the forth and fifth equations is
equivalent to
$$
\partial \psi^2_2 - \bar{\partial} \psi_1^2 +i\psi_1 \psi_2
(|\psi_2|^2 - |\psi_1|^2) = -2H \psi_1 \psi_2 (|\psi_1|^2+|\psi_2|^2).
$$
The system of these two equations at the points where $\psi_1 \psi_2 \neq 0$
are written in the form of the Dirac equation which is satisfied by continuity
on the whole surface:
\begin{equation}
\label{nil-dirac}
\begin{split}
{\cal D}_\nil \psi =
\left[\left(
\begin{array}{cc}
0 & \partial \\
-\bar{\partial} & 0
\end{array}
\right)
 +
\left(
\begin{array}{cc}
U_\nil & 0 \\
0 & V_\nil
\end{array}
\right)
\right]
\psi = 0, \\
U_\nil = V_\nil = \frac{H}{2}(|\psi_1|^2+|\psi_2|^2) +
\frac{i}{4}(|\psi_2|^2-|\psi_1|^2),
\end{split}
\end{equation}
where $H$ is the mean curvature of the surface.

{\sc Remark 1.}
The tangent plane at a point where $Z_3 = \psi_1 \bar{\psi}_2 =0$ is spanned by
the left invariant fields generated by $e_1$ and $e_2$. Since the commutator of this fields
is transversal to their span (i.e., $[e_1,e_2]=e_3$) the left invariant
distribution of $2$-planes
generated by $e_1$ and $e_2$ is nowhere integrable. Therefore
the equality $Z_3=0$ can not hold in open subsets of the surface. We indicate that
the potentials of the Dirac operator are correctly defined everywhere and therefore
the equation ${\cal D}\psi=0$ is extended by continuity to the closure of
$\{Z_3 \neq 0\}$ which, as we argue above,
coincides with the whole surface. The same is true for surfaces in $SU(2)$ and
$\sl$ but fails in the case of $G = \sol$.

The Hopf differential takes the form
\begin{equation}
\label{nil-hopf}
A =  (\bar{\psi}_2 \partial \psi_1 - \psi_1 \partial \bar{\psi}_2) +
i \psi_1^2 \bar{\psi}_2^2
\end{equation}
and the Weingarten equations \eqref{weingarten}
consists of the Dirac equation and
the following system
\begin{equation}
\label{nil-weingarten}
\begin{split}
\partial \psi_1 = \alpha_z \psi_1 + Ae^{-\alpha} \psi_2 -
\frac{i}{2} \psi_1^2 \bar{\psi}_2,
\\
\bar{\partial} \psi_2 = -\bar{A} e^{-\alpha} \psi_1  +
\alpha_{\bar{z}}\psi_2-
 \frac{i}{2} \bar{\psi}_1 \psi_2^2.
\end{split}
\end{equation}

It follows from the Weingarten equations that
$$
(\partial - {\cal A})(\bar{\partial} - {\cal B}) \psi -
(\bar{\partial} - {\cal B})(\partial - {\cal A}) \psi =
({\cal A}_{\bar{z}} - {\cal B}_z + [{\cal A},{\cal B}]) \psi = 0.
$$

{\sc Remark 2.} In the cases $G = \R^3$ or $SU(2)$ the spinor
$\psi^\ast = (-\bar{\psi}_2, \psi_1)^\perp$ meets the equation
$$
R\psi^\ast = ({\cal A}_{\bar{z}} - {\cal B}_z + [{\cal A},{\cal B}])
\psi^\ast = 0
$$
which together with $R\psi=0$ implies that
$$
R = {\cal A}_{\bar{z}} - {\cal B}_z + [{\cal A},{\cal B}] = 0
$$
(see \cite{T3}). In the case of $G=\nil$ as well as for $G=\sl$ or
$\sol$ the equation $R\psi^\ast=0$ does not hold and, in
particular, the kernel of the Dirac operator can not be treated as
a vector space over quaternions. Therefore we shall derive the
Codazzi equations by other methods.

We have
$$
{\cal A} =
\left(
\begin{array}{cc}
\alpha_z - \frac{i}{2}Z_3 & Ae^{-\alpha} \\
-W &  -\frac{i}{2}Z_3
\end{array} \right), \ \ \
{\cal B} =
\left(
\begin{array}{cc}
-\frac{i}{2} \bar{Z}_3 & \overline{W} \\
-\bar{A} e^{-\alpha} &  \alpha_{\bar{z}} - \frac{i}{2} \bar{Z}_3
\end{array} \right)
$$
where $Z_3 = \psi_1 \bar{\psi}_2$ and
$W = \left(\frac{H}{2}-\frac{i}{4} \right)e^\alpha$.
The equation $R\psi = 0$ is written as the system of two equations:
$$
\kappa_1 = ((\alpha_{z \bar{z}} - |A|^2 e^{-2\alpha} + |W|^2)
+ \frac {i}{2}(\partial \bar{Z}_3 - \bar{\partial}
Z_3))\psi_1 +
(A_{\bar{z}} e^{-\alpha} - \overline{W}_z  +
\alpha_z \overline{W}) \psi_2 = 0,
$$
$$
\kappa_2 = (\bar{A}_z e^{-\alpha} - W_{\bar{z}}  +
\alpha_{\bar{z}}W ) \psi_1 +
(-(\alpha_{z \bar{z}} - |A|^2 e^{-2\alpha} + |W|^2)
+ \frac {i}{2}(\partial \bar{Z}_3 - \bar{\partial}
Z_3))\psi_2 = 0.
$$
We recall that surfaces in $\nil$ we have
\begin{equation}
\label{nil-identity}
\partial \bar{Z}_3 - \bar{\partial} Z_3 =
- \frac{i}{2} (|\psi_2|^4-|\psi_1|^4), \ \ \
\partial \bar{Z}_3 + \bar{\partial}  Z_3=
H (|\psi_2|^4 - |\psi_1|^4).
\end{equation}
The Codazzi equations is now the system
$$
\kappa_1 \bar{\psi}_1 - \bar{\kappa}_2 \psi_2 = 0, \ \ \
\kappa_1 \bar{\psi}_2 + \bar{\kappa}_2 \psi_1 = 0
$$
which is rewritten as
$$
\alpha_{z\bar{z}} - |A|^2 e^{-2\alpha} + \frac{H^2}{4}
e^{2\alpha} = \frac{1}{16}(3|\psi_1|^4 + 3 |\psi_2|^4 - 10 |\psi_1|^2
|\psi_2|^2),
$$
$$
A_{\bar{z}} - \frac{H_z}{2}e^{2\alpha} +
\frac{1}{2}(|\psi_2|^4-|\psi_1|^4)\psi_1 \bar{\psi}_2 = 0.
$$
By using \eqref{nil-identity} we rewrite this system as follows
\begin{equation}
\label{nil-codazzi}
\begin{split}
\alpha_{z \bar{z}} -  e^{-2 \alpha} |A|^{2} +
\frac{1}{4} e^{2\alpha} H^{2} =
\frac{3}{16}e^{2\alpha} - |Z_3|^2, \\
\bar{\partial}\left (A + \frac{{Z_3}^2}{2H+i} \right)
= \frac{1}{2} H_z e^{2\alpha} +
\bar{\partial} \left(\frac{1}{2H+i} \right) {Z_3}^2.
\end{split}
\end{equation}

We summarize the above results in the following

\begin{theorem}
For a surface in $G=\nil$ its generating spinor
$\psi$ meets the Dirac equation \eqref{nil-dirac}.

Moreover each function $\psi$ meeting \eqref{nil-dirac}
serves as the generating spinor of a surface in $\nil$.

The Weingarten equations for the surface are given by \eqref{nil-dirac}
and \eqref{nil-weingarten}.
The Hopf differential takes the form \eqref{nil-hopf}
and the Codazzi equations are given by \eqref{nil-codazzi}.
\end{theorem}

\begin{corollary}
The generating spinor $\psi$
of a minimal surface in $\nil$ meets the following
equations
$$
\bar{\partial} \psi_1 = \frac{i}{4}(|\psi_2|^2-|\psi_1|^2)\psi_1, \ \ \
\partial \psi_2 = - \frac{i}{4}(|\psi_2|^2-|\psi_1|^2) \psi_1.
$$
\end{corollary}

\begin{corollary}[Abresch]
For a constant mean curvature surface in $\nil$ the quadratic
differential
$$
\widetilde{A} dz^2 = \left (A + \frac{{Z_3}^2}{2H+i} \right)dz^2
$$
is holomorphic.
\end{corollary}

\begin{proposition}
\label{nil-holomorphic}
If the differential $\widetilde{A} dz^2$ is holomorphic then the surface in $\nil$
has constant mean curvature.
\end{proposition}

{\sc Proof.}
We prove that by contradiction. Let us assume that $H_z \neq 0$ in some domain.
Then we have
$$
\frac{1}{2}e^{2\alpha} = \frac{2H_{\bar{z}}}{(2H+i)^2} Z_3^2
$$
and which implies the equality for the moduli of both sides:
$$
\frac{1}{2}e^{2\alpha} = \frac{2 |\psi_1|^2 |\psi_2|^2}{4H^2+1}.
$$
Since $e^\alpha = |\psi_1|^2 + |\psi_2|^2$, this is rewritten as
$$
(4H^2+1)(|\psi_1|^+|\psi_2|^2)^2 = 4|\psi_1|^2|\psi_2|^2
$$
which implies
$$
4H^2 e^{2\alpha} + (|\psi_1|^2 - |\psi_2|^2)^2 = 0.
$$
It is clear that this equality is valid if and only if $|\psi_1| = |\psi_2|$ and $H = 0$,
i.e. the surface is minimal. This contradiction proves the proposition.

Finally we compute the energy of a compact surface in $\nil$.

\begin{proposition}
\label{nil-energy-prop}
Given a closed oriented surface $M$ in the group $\nil$,
its energy is real-valued and equals
\begin{equation}
\label{nil-energy}
E(M) = \int_M \left( \frac{H^2}{4}(|\psi_1|^2+
|\psi_2|^2)^2 - \frac{1}{16}(|\psi_2|^2-|\psi_1|^2)^2\right)
\frac{idz\wedge d\bar{z}}{2}
\end{equation}
\end{proposition}

{\sc Proof.} Let us insert the formulas for $U_\nil$ and $V_\nil$ into
\eqref{eq-identity} and obtain
$$
\bar{\partial} (\psi_1\bar{\psi}_2) =
\frac{H}{2}(|\psi_2|^2 - |\psi_1|^2) + i \frac{1}{4}(|\psi_2|^4-|\psi_1|^4).
$$
By the Stokes formula we have
\begin{equation}
\label{nil-equality-0}
\int_M \frac{H}{2}(|\psi_2|^2 - |\psi_1|^2)
\frac{idz\wedge d\bar{z}}{2} + i \int_M \frac{1}{4} (|\psi_2|^4 -
|\psi_1|^4) \frac{idz\wedge d\bar{z}}{2} = 0.
\end{equation}
The real part of the left hand side of this formula is a multiple of
$\Im E(M)$ which implies that the energy is real-valued.
We insert now $U_\nil$ and $V_\nil$
into $\Re E(M)$ and derive \eqref{nil-energy}. This proves the proposition.

Notice that the vanishing of the the imaginary part of the left hand side of
\eqref{nil-equality-0} together with \eqref{normal} implies
the equality
$$
\int_M \langle f^{-1}(N), e_3 \rangle d\mu =
\int_M (|\psi_2|^4 - |\psi_1|^4)\frac{idz\wedge d\bar{z}}{2} = 0.
$$

\begin{proposition}
The energy of a surface $M$ in the group $\nil$ equals
\begin{equation}
\label{nil-energy-geo}
E(M) = \frac{1}{4} \int_M \left( H^2 + \frac{\widehat{K}}{4} - \frac{1}{16} \right)
d \mu,
\end{equation}
where $\widehat{K}$ is the sectional curvature of the tangent plane at a point.
\end{proposition}

{\sc Proof.} It is computed by using Proposition \ref{curvatures} that the sectional curvature of a tangent plane
at the unit of group equals
$$
\widehat{K} = \frac{1}{4} - \cos^2 \varphi
$$
where $\varphi$ is the angle between the normal to the plane and $e_3$.
By \eqref{normal}, the integrand in \eqref{nil-energy} is equal to
$$
\frac{1}{4}e^{2\alpha} H^2 - \frac{1}{16}(|\psi_2|^2 - |\psi_1|^2)^2 =
\frac{1}{4}e^{2\alpha} H^2 - \frac{1}{16}e^{2\alpha} \langle f^{-1}(N),e_3\rangle^2 =
$$
$$
\frac{1}{4}e^{2\alpha} \left(H^2 - \frac{1}{4} \cos^2 \varphi \right) =
\frac{1}{4}e^{2\alpha} \left(H^2 - \frac{1}{4} \left(\frac{1}{4} - \widehat{K}\right)\right) =
$$
$$
= \frac{1}{4} e^{2\alpha} \left(H^2 + \frac{1}{4} \widehat{K} -\frac{1}{16}\right).
$$
The proposition is proved.

\subsection{The group $\sl$}
\label{ssec3.2}

Inserting \eqref{sl-con} into  the derivational equations \eqref{z1} and
\eqref{z2} we obtain
\begin{equation}
\label{sl-1}
\begin{aligned}
\partial \bar{Z}_1 - \bar{\partial} Z_1 +
\frac{1}{2}(Z_2 \bar{Z}_3 - \bar{Z}_2 Z_3) = 0, \\
\partial \bar{Z}_2 - \bar{\partial} Z_2 +
\frac{1}{2}(Z_3 \bar{Z}_1 - \bar{Z}_3 Z_1) = 0, \\
\partial \bar{Z}_3 - \bar{\partial} Z_3 -
2(Z_1\bar{Z}_2 - \bar{Z}_1 Z_2) = 0,\\
\partial \bar{Z}_1 + \bar{\partial} Z_1 -
\frac{5}{2}(Z_2\bar{Z}_3 + \bar{Z}_2 Z_3) =
2iH(\bar{Z}_2 Z_3 - Z_2 \bar{Z}_3), \\
\partial \bar{Z}_2 + \bar{\partial} Z_2 +
\frac{5}{2}(Z_1\bar{Z}_3 + \bar{Z}_1 Z_3)=
2iH(\bar{Z}_3 Z_1 - Z_3 \bar{Z}_1), \\
\partial \bar{Z}_3 + \bar{\partial} Z_3 =
2iH(\bar{Z}_1 Z_2 - Z_1 \bar{Z}_2).
\end{aligned}
\end{equation}
By inserting \eqref{spinor} into these equations we reduce the system of
the first two of them to
$$
\partial \psi_2^2 + \bar{\partial} \psi_1^2 -\frac{i}{2}\psi_1 \psi_2
(|\psi_1|^2+|\psi_2|^2) = 0,
$$
and the fourth and fifth equations together are reduced to
$$
\partial \psi_2^2 - \bar{\partial} \psi_1^2 +\frac{5}{2}i\psi_1\psi_2
(|\psi_1|^2-|\psi_2|^2) = -2H \psi_1 \psi_2 (|\psi_1|^2+|\psi_2|^2).
$$
These two equations are written in terms of
the Dirac equation (the derivation is completely similar to the $\nil$
case, see \eqref{nil-dirac} in \S \ref{ssec3.1}]:
\begin{equation}
\label{sl-dirac}
\begin{split}
{\cal D}_\sll \psi =
\left[\left(
\begin{array}{cc}
0 & \partial \\
-\bar{\partial} & 0
\end{array}
\right)
 +
\left(
\begin{array}{cc}
U_\sll & 0 \\
0 & V_\sll
\end{array}
\right)
\right]
\psi = 0, \\
U_\sll = \frac{H}{2}(|\psi_1|^2+|\psi_2|^2) +
i\left(\frac{1}{2}|\psi_1|^2-\frac{3}{4}|\psi_2|^2\right),
\\
V_\sll = \frac{H}{2}(|\psi_1|^2+|\psi_2|^2) +
i\left(\frac{3}{4}|\psi_1|^2-\frac{1}{2}|\psi_2|^2\right),
\end{split}
\end{equation}
where $H$ is the mean curvature of the surface.

The Hopf differential equals
\begin{equation}
\label{sl-hopf}
A =  (\bar{\psi}_2 \partial \psi_1 - \psi_1 \partial \bar{\psi}_2)
-\frac{5i}{2}  \psi_1^2 \bar{\psi}_2^2
\end{equation}
and we have
\begin{equation}
\label{sl-weingarten}
\begin{split}
\partial \psi_1 = \alpha_z \psi_1 + Ae^{-\alpha} \psi_2 +
\frac{5i}{4} \psi_1^2 \bar{\psi}_2,
\\
\bar{\partial} \psi_2 = -\bar{A} e^{-\alpha} \psi_1  +
\alpha_{\bar{z}}\psi_2 +
 \frac{5i}{4} \bar{\psi}_1 \psi_2^2.
\end{split}
\end{equation}
The matrices ${\cal A}$ and ${\cal B}$ are given by the formulas
$$
{\cal A} =
\left(
\begin{array}{cc}
\alpha_z + \frac{5i}{4}Z_3 & Ae^{-\alpha} \\
-W &  \frac{5i}{4}Z_3
\end{array} \right), \ \ \
{\cal B} =
\left(
\begin{array}{cc}
\frac{5i}{4} \bar{Z}_3 & \overline{W} \\
-\bar{A} e^{-\alpha} &  \alpha_{\bar{z}} + \frac{5i}{4} \bar{Z}_3
\end{array} \right)
$$
where $W = \frac{1}{2}(H+i)e^\alpha$.
The equation $R\psi = 0$ is equivalent to the equations
$$
\kappa_1 = ((\alpha_{z \bar{z}} - |A|^2 e^{-2\alpha} + |W|^2)
- \frac {5i}{4}(\partial \bar{Z}_3 - \bar{\partial}
Z_3))\psi_1 +
(A_{\bar{z}} e^{-\alpha} - \overline{W}_z  +
\alpha_z \overline{W}) \psi_2 = 0,
$$
$$
\kappa_2 = (\bar{A}_z e^{-\alpha} - W_{\bar{z}}  +
\alpha_{\bar{z}}W ) \psi_1 +
(-(\alpha_{z \bar{z}} - |A|^2 e^{-2\alpha} + |W|^2)
- \frac {5i}{4}(\partial \bar{Z}_3 - \bar{\partial}
Z_3))\psi_2 = 0.
$$
Since the generating spinors of surfaces in $\sl$
meet the equalities
$$
\partial \bar{Z}_3 - \bar{\partial} Z_3 =
i(|\psi_2|^4-|\psi_1|^4), \ \ \
\partial \bar{Z}_3 + \bar{\partial}  Z_3=
H (|\psi_2|^4 - |\psi_1|^4).
$$
we rewrite the equations
$\kappa_1 \bar{\psi}_1 - \bar{\kappa}_2 \psi_2 = 0$ and
$\kappa_1 \bar{\psi}_2 + \bar{\kappa}_2 \psi_1 = 0$ as the following system
\begin{equation}
\label{sl-codazzi}
\begin{split}
\alpha_{z \bar{z}} -  e^{-2 \alpha} |A|^{2} +
\frac{1}{4} e^{2\alpha} H^{2} =
e^{2\alpha} - 5|Z_3|^2, \\
\bar{\partial}\left (A + \frac{5 Z_3^2}{2(H-i)} \right)
= \frac{1}{2} H_z e^{2\alpha}
+ \bar{\partial} \left(\frac{5}{2(H-i)} \right) {Z_3}^2.
\end{split}
\end{equation}

We derive

\begin{theorem}
Given a surface in $G=\sl$, its generating spinor
$\psi$ satisfies the Dirac equation \eqref{sl-dirac}.

Any function $\psi$ meeting \eqref{sl-dirac}
is the generating spinor of a surface in $\sl$.

The Weingarten equations for the surface is the system consisting of
\eqref{sl-dirac}
and \eqref{sl-weingarten}.
The Hopf differential is given by \eqref{nil-hopf}
and the Codazzi equations take the form \eqref{nil-codazzi}.
\end{theorem}

\begin{corollary}
The generating spinor of a minimal surface in $\sl$ satisfies
the equations
$$
\bar{\partial} \psi_1 =
i\left(\frac{3}{4}|\psi_1|^2-\frac{1}{2}|\psi_2|^2\right)
\psi_2,
\ \ \
\partial \psi_2 =
-i\left(\frac{1}{2}|\psi_1|^2-\frac{3}{4}|\psi_2|^2\right)\psi_1.
$$
\end{corollary}

\begin{corollary}
[Abresch]
For a constant mean curvature surface in $\sl$ the quadratic differential
$$
\widetilde{A} dz^2 = \left(A + \frac{5}{2(H-i)}Z^2_3\right) dz^2
$$
is holomorphic.
\end{corollary}

{\sc Remark 3.} The method of proving Proposition
\ref{nil-holomorphic} fails in the case of $\sl$ and we do
not know whether the holomorphicity of $\widetilde{A}$ implies that the
surface has constant mean curvature.

We finish our study of surfaces in $\sl$ by computing the energy functional.

\begin{proposition}
\label{sl-energy-prop}
Given a closed oriented surface $M$ in the group $\sl$,
its energy is real-valued and equals
\begin{equation}
\label{sl-energy}
\begin{split}
E(M) = \int_M \left[ \frac{H^2}{4}(|\psi_1|^2+ |\psi_2|^2)^2 -
\right.
\\
\left.
\left(\frac{1}{2}|\psi_1|^2-\frac{3}{4}|\psi_2|^2\right)
\left(\frac{3}{4}|\psi_1|^2-\frac{1}{2}|\psi_2|^2\right)
\right]
\frac{idz\wedge d\bar{z}}{2}.
\end{split}
\end{equation}
\end{proposition}

The proof of this proposition is straightforward and analogous to the proof of
Proposition \ref{nil-energy-prop}. Hence we skip it and only mention that
as for surfaces in $\nil$ we have
$\int_M \langle f^{-1}(N), e_3 \rangle d\mu =
\int_M (|\psi_2|^4 - |\psi_1|^4)\frac{idz\wedge d\bar{z}}{2} = 0$.

As for the group $\nil$ we can write down the formula for the energy
in common geometric terms.

\begin{proposition}
The energy of a surface $M \subset \sl$ equals
\begin{equation}
\label{sl-energy-geo}
E(M) = \frac{1}{4} \int_M \left( H^2 +
\frac{5}{16} \widehat{K} - \frac{1}{4} \right) d \mu,
\end{equation}
where $\widehat{K}$ is the sectional curvature of the tangent plane at a point.
\end{proposition}

This proposition is straightforward from \eqref{nil-energy} and Proposition
\ref{curvatures}.

\subsection{The group $\sol$}
\label{ssec3.3}

By \eqref{sol-con}, the derivational equations \eqref{z1} and
\eqref{z2} for surfaces in the group $\sol$ take the form
\begin{equation}
\label{sol-1}
\begin{split}
\partial \bar{Z}_1 - \bar{\partial} Z_1 +
(Z_1 \bar{Z}_3 - \bar{Z}_1 Z_3) = 0, \\
\partial \bar{Z}_2 - \bar{\partial} Z_2 -
(Z_2 \bar{Z}_3 - \bar{Z}_2 Z_3) = 0, \\
\partial \bar{Z}_3 - \bar{\partial} Z_3 = 0,\\
\partial \bar{Z}_1 + \bar{\partial} Z_1 + (Z_1\bar{Z}_3 + \bar{Z}_1 Z_3) =
2iH(\bar{Z}_2 Z_3 - Z_2 \bar{Z}_3), \\
\partial \bar{Z}_2 + \bar{\partial} Z_2 - (Z_2\bar{Z}_3 + \bar{Z}_2 Z_3)=
2iH(\bar{Z}_3 Z_1 - Z_3 \bar{Z}_1), \\
\partial \bar{Z}_3 + \bar{\partial} Z_3 -2(|Z_1|^2-|Z_2|^2) =
2iH(\bar{Z}_1 Z_2 - Z_1 \bar{Z}_2).
\end{split}
\end{equation}

As in \S \ref{ssec3.1} and \S \ref{ssec3.2}
we insert \eqref{spinor} in these equations and
derive the following pair of equations
$$
\partial \psi_2^2 + \bar{\partial} \psi_1^2 -\bar{\psi}_1 \bar{\psi}_2
(|\psi_1|^2+|\psi_2|^2) = 0,
$$
$$
\partial \psi_2^2 - \bar{\partial} \psi_1^2 + \bar{\psi}_1\bar{\psi}_2
(|\psi_1|^2-|\psi_2|^2) = -2H \psi_1 \psi_2 (|\psi_1|^2+|\psi_2|^2).
$$

We see that at points where $Z_3 = \psi_1 \bar{\psi}_2 \neq 0$ these equations are
rewritten in terms of the Dirac equation
\begin{equation}
\label{sol-dirac}
\begin{split}
{\cal D}_\sol \psi =
\left[\left(
\begin{array}{cc}
0 & \partial \\
-\bar{\partial} & 0
\end{array}
\right)
 +
\left(
\begin{array}{cc}
U_\sol & 0 \\
0 & V_\sol
\end{array}
\right)
\right]
\psi = 0, \\
U_\sol = \frac{H}{2}(|\psi_1|^2+|\psi_2|^2) +
\frac{1}{2}\bar{\psi}_2^2 \frac{\bar{\psi}_1}{\psi_1},
\\
V_\sol = \frac{H}{2}(|\psi_1|^2+|\psi_2|^2) +
\frac{1}{2}\bar{\psi}_1^2 \frac{\bar{\psi}_2}{\psi_2}.
\end{split}
\end{equation}

{\sc Remark 4.} Since the left invariant vector fields generated
by $e_1$ and $e_2$ commute the equation $Z_3 = \psi_1 \bar{\psi}_2
=0$ can be valid in an open subset $B$ of a surface. Therefore the
Dirac equation can not be extended by continuity onto the whole
surface and does not describe $\psi$ in $B$. Since $H=0$ in $B$ it
is reasonable to assume that $U_\sol = V_\sol=0$ inside such a
domain. However the
potentials $U_\sol$ and $V_\sol$ are  not always correctly defined
on the boundaries of the set $\{Z_3 \neq 0\}$
in view of the indeterminacy of $\frac{\bar{\psi}_1}{\psi_1}$ and
$\frac{\bar{\psi}_2}{\psi_2}$. The measure of the indeterminacy
set is zero and we may correctly define the energy of a surface in
$\sol$ as
$$
E(M) = \int_{\{Z_3 \neq 0\}} U_\sol V_\sol d\mu.
$$
We did not manage to rewrite this quantity in common geometric terms as we did for surfaces
in $\nil$ and $\sl$ and its geometric meaning remains still unclear for us.
We even do not know whether the energy real-valued (as for surfaces in $SU(2), \nil$, and $\sl$)
or not.

The Hopf differential of a surface in $\sol$ equals
\begin{equation}
\label{sol-hopf}
A =  (\bar{\psi}_2 \partial \psi_1 - \psi_1 \partial \bar{\psi}_2) +
\frac{1}{2} (\bar{\psi}_2^4 - \psi_1 ^4)
\end{equation}
and we complete the Weingarten equations by the following system
\begin{equation}
\label{sol-weingarten}
\begin{split}
\partial \psi_1 = \alpha_z \psi_1 + Ae^{-\alpha} \psi_2 -
\frac{1}{2} \bar{\psi}_2^3, \\
\bar{\partial} \psi_2 = -\bar{A} e^{-\alpha} \psi_1  + \alpha_{\bar{z}}\psi_2-
\frac{1}{2} \bar{\psi}_1^3.
\end{split}
\end{equation}
The Codazzi equations take the form
\begin{equation}
\label{sol-codazzi}
\begin{split}
\alpha_{z \bar{z}} -  e^{-2 \alpha} |A|^{2} +
\frac{1}{4} e^{2\alpha} H^{2} = \frac{1}{4}
(6 |\psi_1|^2 |\psi_2|^2  -(|\psi_1|^4 +  |\psi_2|^4)), \\
A_{\bar{z}} - \frac{1}{2} H_z e^{2\alpha} =
(|\psi_2|^4 - |\psi_1|^4) \psi_1 \bar{\psi}_2.
\end{split}
\end{equation}

Let us summarize these results.
Let $f: M \to \sol$ be a surface. We
by $B$ the subset of $M$ where $Z_3 = \langle f^{-1}f_z, e_3
\rangle =0$, by $B_0$ the interior of $B$, and by $C$ the subset of $M$
where $Z_3 \neq 0$. Then $M = B_0 \cup \bar{C}$ and the set
$B \setminus B_0$ lies in the closure $\bar{C}$ of $C$ and has
measure zero.

\begin{theorem}
The generating spinor $\psi$ of a surface $f: M \to \sol$ meets the Dirac equation
\eqref{sol-dirac} in $C$ and meets the Dirac equation with zero potentials:
$\bar{\partial}\psi_1 = \partial \psi_2 =0$, in $B_0$.

Every function $\psi$ meeting \eqref{sol-dirac} is some set $D \subset M$ is
the generating function of a surface $f: D \to \sol$.

The Hopf differential is given by \eqref{sol-hopf}.
In $B_0 \cup C$ the Weingarten equations are formed by \eqref{sol-weingarten} and by
the Dirac equation for $\psi$.
The Codazzi equations take the form \eqref{sol-codazzi}.
\end{theorem}

\begin{corollary}
The generating spinor $\psi$ of a minimal surface in $\sol$
meets the equations
$$
\bar{\partial}\psi_1 = \frac{1}{2} \bar{\psi}_1^2 \bar{\psi}_2, \ \ \
\partial \psi_2 = - \frac{1}{2} \bar{\psi}_1 \bar{\psi}_2^2.
$$
\end{corollary}

\medskip

{\bf Addition to the proofs.} After submitting this paper to the
journal we became aware of the papers \cite{I1,I2}, where the
authors obtained some interesting results on minimal surfaces in
three-dimensional Lie groups under study. In particular, some analogs
of the Weierstrass representation for minimal surfaces in $\nil$
and  $\sol$ were derived.

\end{document}